\documentclass[12pt]{article}
\usepackage{amssymb}
\usepackage{amsfonts}
\usepackage{latexsym}
\usepackage[dvips]{graphics}

\large\normalsize

\begin{document}
\begin{center}
{\Large\bf Enumeration of perfect matchings of a type of Cartesian
products of graphs$^\ast$}
\\[30pt]
{Weigen\ Yan$^{\rm a,b}$ \footnote{This work is supported by
FMSTF(2004J024) and FJCEF(JA03131)} \quad and \quad Fuji\
Zhang$^{\rm b}$ \footnote{Partially supported by NSFC(10371102).
\newline\hspace*{5mm}{\it Email address:} weigenyan@263.net (W.
Yan), fjzhang@jingxian.xmu.edu.cn (F. Zhang)}}
\\[10pt]
\footnotesize { $^{\rm a}$School of Sciences, Jimei University,
Xiamen 361021, China
\\[7pt]
$^{\rm b}$Department of Mathematics, Xiamen University, Xiamen
361005, China}
\end{center}
\begin{abstract}
Let $G$ be a graph and let Pm$(G)$ denote the number of perfect
matchings of $G$. We denote the path with $m$ vertices by $P_m$
and the Cartesian product of graphs $G$ and $H$ by $G\times H$. In
this paper, as the continuance of our paper [19], we enumerate
perfect matchings in a type of Cartesian products of graphs by the
Pfaffian method, which was discovered by Kasteleyn.
Here are some of our results:\\
1. Let $T$ be a tree and let $C_n$ denote the cycle with $n$
vertices. Then Pm$(C_4\times T)=\prod (2+\alpha^2)$, where the
product ranges over all eigenvalues $\alpha$ of $T$.
Moreover, we prove that Pm$(C_4\times T)$ is always a square or double a square.\\
2. Let $T$ be a tree. Then Pm$(P_4\times T)=\prod
(1+3\alpha^2+\alpha^4)$, where the product ranges over all
non-negative
eigenvalues $\alpha$ of $T$.\\
3. Let $T$ be a tree with a perfect matching. Then Pm$(P_3\times
T)=\prod (2+\alpha^2),$ where the product ranges over all positive
eigenvalues $\alpha$ of $T$. Moreover, we prove that Pm$(C_4\times
T)=[\mbox{Pm}(P_3\times T)]^2$.
\end{abstract}
\section* { \bf 1. Introduction }
\par A {\bf perfect matching} of a graph $G$ is a set of independent edges of $G$ covering all vertices
of $G$. Problems involving enumeration of perfect matchings of a graph were first examined by chemists
and physicists in the 1930s (for history see [4,16]), for two different (and unrelated) purposes: the study
of aromatic hydrocarbons and the attempt to create a theory of the liquid state.
\par Shortly after the advent of quantum chemistry, chemists turned their attention to molecules like
benzene composed of carbon rings with attached hydrogen atoms. For
these researchers, perfect matchings of a polyhex graph
corresponded to "Kekul\'e structures", i.e., assigning single and
double bonds in the associated hydrocarbon (with carbon atoms at
the vertices and tacit hydrogen atoms attached to carbon atoms
with only two neighboring carbon atoms). There are strong
connections between combinatorial and chemical properties for such
molecules; for instance, those edges which are present in
comparatively few of the perfect matchings of a graph turn out to
correspond to the bonds that are least stable, and the more
perfect matchings a polyhex graph possesses the more stable is the
corresponding benzenoid molecule. The number of perfect matchings
is an important topological index which had been applied for
estimation of the resonant energy and total $\pi-$electron energy
and calculation of pauling bond order (see [6,15,18]). So far,
many mathematicians, physicists and chemists have given most of
their attention to counting perfect matchings of graphs. See for
example papers [2,5,7,16,17,19$-$23].
\par By a simple graph $G=(V(G),E(G))$ we mean a finite undirected graph, that is, one with no loops or
parallel edges, with the vertex$-$set
$V(G)=\{v_1,v_2,\ldots,v_n\}$ and the edge$-$set
$E(G)=\{e_1,e_2,\ldots, e_m\}$, if not specified. We denote by
${\bf \mbox{Pm}(G)}$ the number of perfect matchings of $G$. If
$M$ is a perfect matching of $G$, an {\bf $ M-$alternating cycle}
in $G$ is a cycle whose edges are alternately in $E(G)\backslash
M$ and $M$. Let $G$ be a graph. A cycle $C$ of $G$ is called to be
{\bf nice} if $G-C$ contains a perfect matching, where $G-C$
denotes the induced subgraph of $G$ obtained from $G$ by deleting
the vertices of $C$. Throughout this paper, we denote a tree by
$T$ and a path with $n$ vertices by $P_n$. For two graphs $G$ and
$H$, let $G\times H$ denote the Cartesian product of graphs $G$
and $H$.
\par Let $G=(V(G),E(G))$ be a simple graph and let $G^e$ be
an arbitrary orientation of $G$.  The {\bf skew adjacency matrix} of $G^e$, denoted by ${\bf A(G^e)}$, is defined
as follows:
$$
A(G^e)=(a_{ij})_{n\times n},
$$
$$
a_{ij}=\left\{\begin{array}{ll}
1 & \mbox{if} \ \  (v_i,v_j)\in E(G^e),\\
-1 & \mbox{if} \ \  (v_j,v_i)\in E(G^e),\\
0 &   \mbox{otherwise}.
\end{array}
\right.
$$
It is clear that
$A(G^e)$ is a skew symmetric matrix, that is, $A(G^e)^T=-A(G^e)$.
\par Let $G$ be a simple graph. We say $G$ has {\bf reflective symmetry} if it is invariant under the reflection
across some straight line or plane $l$ (the symmetry plane or
axis) (see Ciucu's paper [2]). Ciucu [2] gave a matching
factorization theorem of the number of perfect matchings of a
symmetric plane bipartite graph in which there are some vertices
lying on the symmetry axis $l$ but no edges crossing $l$. Ciucu's
theorem expresses the number of perfect matchings of $G$ in terms
of the product of the number of perfect matchings of two subgraphs
of $G$ each one of which has nearly half the number of vertices of
$G$. On the other hand, in [21] Zhang and Yan proved that if a
bipartite graph $G$ without nice cycles of length $4s, s\in
\{1,2,\ldots\}$ was invariant under the reflection across some
plane (or straight line) and there are no vertices lying on the
symmetry plane (or axis) then $\mbox{Pm}(G)=|\det A(G^+)| $, where
$G^+$ is a graph having loops with half the number of vertices of
$G$, and $A(G^+)$ is the adjacency matrix
of $G^+$. Furthermore, in [19] Yan and Zhang obtained the following results on the symmetric graphs:\\
1. If $G$ is a reflective symmetric plane graph (which does not need to be bipartite)
without vertices on the symmetry axis, then the number of perfect
matchings of $G$ can be expressed by a determinant of order $\frac {1}{2}|G| $, where $|G|$ denotes
the number of vertices of $G$.\\
2. Let $G$ be a bipartite graph without cycles of length $4s, s\in
\{1, 2, \ldots \}$. Then the number of perfect matchings of
$G\times P_2$ equals $\prod (1+\alpha^2)$, where the product
ranges over all non-negative eigenvalues $\alpha$ of $G$ .
Particularly, if $T$ is a tree then $\mbox{Pm}(T\times P_2)$
equals $\prod (1+\theta^2)$, where the product ranges over all non-negative eigenvalues $\theta$ of $T$. \\
\par As the continuance of our paper [19],
in this paper we obtain the following results:\\
1. Let $T$ be a tree and let $C_n$ denote the cycle with $n$
vertices. Then $\mbox{Pm}(C_4\times T)=\prod (2+\alpha^2)$, where
the product ranges over all of eigenvalues $\alpha$ of $T$. This
makes it possible to obtain a formula for the number of perfect
matchings for the linear $2\times 2\times n$ cubic lattice, which
was previously obtained by H. Narumi and H. Hosoya [14]. Moreover,
we prove that $\mbox{Pm}(C_4\times T)$ is always a square or
double a square (such a number
was called squarish in [7]).\\
2. Let $T$ be a tree. Then $\mbox{Pm}(P_4\times T)=\prod
(1+3\alpha^2+\alpha^4)$, where the product ranges over all
non-negative
eigenvalues $\alpha$ of $T$.\\
3. Let $T$ be a tree with a perfect matching. Then
$\mbox{Pm}(P_3\times T)=\prod (2+\alpha^2),$ where the product
ranges over all positive eigenvalues $\alpha$ of $T$. Moreover, we
prove that
$\mbox{Pm}(C_4\times T)=[\mbox{Pm}(P_3\times T)]^2$.\\
\par The start point of this paper is the fact that we can use the Pfaffian method,
which was discovered by Kasteleyn [8,9,11], to enumerate perfect
matchings of some graphs. In order to formulate lemmas we need to
introduce some terminology and notation as follows.
\par If $D$ is an orientation of a simple graph $G$ and $C$ is a cycle of even length, we say that $C$ is {\bf oddly
oriented} in $D$ if $C$ contains odd number of edges that are
directed in $D$ in the direction of each orientation of $C$ (see
[5,11]). We say that $D$ is a {\bf Pfaffian orientation} of $G$ if
every nice cycle of even length of $G$ is oddly oriented in $D$.
It is well known that if a graph $G$ contains no subdivision of
$K_{3,3}$ then $G$ has a Pfaffian orientation (see Little [10]).
McCuaig [12], and McCuaig, Robertson et al [13], and Robertson,
Seymour et al [17]
found a polynomial-time algorithm to determine whether a bipartite graph has a Pfaffian orientation.\\
{\bf Lemma 1 [8,9,11]} Let $G^e$ be a Pfaffian orientation of a graph $G$. Then\\
$$
[\mbox{Pm}(G)]^2=\det A(G^e),
$$
where $A(G^e)$ is the skew adjacency matrix of $G^e$.\\
{\bf Lemma 2 [11]} Let $G$ be any simple graph with even number of
vertices, and $G^e$ an orientation of $G$. Then the
following three properties are equivalent:\\
(1) $G^e$ is a Pfaffian orientation.\\
(2) Every nice cycle of even length in $G$ is oddly oriented in $G^e$.\\
(3) If $G$ contains a perfect matching, then for some perfect
matching $F$, every $F$-alternating cycle is oddly oriented in
$G^e$.
\begin{figure}[htbp]
  \centering
  \scalebox{0.75}{\includegraphics{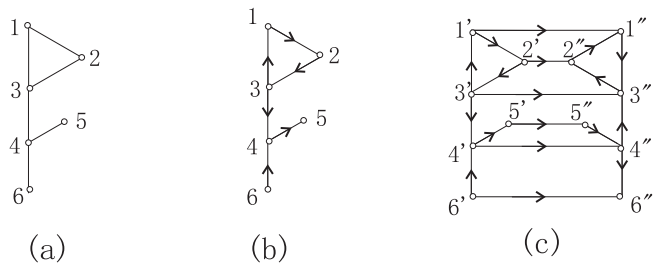}}
  \caption{\ (a)\ A graph $G$.
  \ (b)\ An orientation $G^e$ of $G$. (c)\ The orientation $(P_2\times G)^e$ of $P_2\times G$.}
\end{figure}
\section* { \bf 2. Enumeration of perfect matchings of $C_4\times T$ }
\par First, we introduce a method to orient a type of symmetric graphs. Let $G$ be a simple graph, and $G^e$ an
orientation of $G$. We take a copy of $G^e$, denoted by $G_1^e$.
If we reverse the orientation of each arc of $G_1^e$, then we
obtain another orientation of $G$, denoted by $G_2^e$. Hence
$G_2^e$ is the converse of $G_1^e$. Note that $P_2\times G$ can be
obtained as follows: Take two copies of $G$, denoted by $G_1$ with
vertex-set $V(G_1)=\{v'_1, v'_2, \ldots, v'_n\}$ and $G_2$ with
vertex-set $V(G_2)=\{v''_1, v''_2, \ldots, v''_n\}$ (we consider
that $G_1$ and $G_2$ are {\bf the left half} and {\bf right half}
of $P_2\times G$, respectively), and add an edge $v'_iv''_i$
between every pair of corresponding vertices $v'_i$ and $v''_i$ of
$G_1$ and $G_2$, respectively. It is obvious that the resulting
graph is $P_2\times G$ and all edges $v'_iv''_i$ (for $1\leq i\leq
n$) added between the left half and the right half of  $P_2\times
G$ form a perfect matching of $P_2\times G$, denoted by $M$,  and
$G_1^e $ (or $G_2^e$) is an orientation of $G_1$ (or $G_2$). If we
define the direction of every edge in $M$ is from the left to the
right, then an orientation of $P_2\times G$ is obtained, denoted
by {\bf $(P_2\times G)^e$}.  Figure 1 illustrates
this procedure. By using a result from Fischer and Little [5], the following lemma was proved by Yan and Zhang in [19].\\
{\bf Lemma 3 [19]} Let $G$ be a simple graph. If $G^e$ is an
orientation of $G$ under which every cycle of even length is oddly
oriented in $G^e$, then the orientation $(P_2\times G)^e$ defined
as above is a Pfaffian orientation of $P_2\times G$.
\par Now we prove the following lemma.\\
{\bf Lemma 4} Let $T$ be a tree. Then every cycle of $P_2\times T$ is a nice cycle.\\
{\bf Proof} Let $T_1$ and $T_2$ denote the left half and the right half of $P_2\times T$, respectively. Suppose that
$C$ is a cycle of $P_2\times T$. We claim that $C$ has the following form:\\
$$v_{i_1}'-v_{i_2}'-\ldots-v_{i_m}'-v_{i_m}''-v_{i_{m-1}}''-\ldots-v_{i_2}''-v_{i_1}''-v_{i_1}',$$
where $v_{i_1}', v_{i_2}', \ldots, v_{i_m}'\in V(T_1)$ and
$v_{i_1}'', v_{i_2}'', \ldots, v_{i_m}''\in V(T_2)$.
\par Note that there is a unique path between two vertices of tree $T$.
Hence we only need to prove that $|E(C)\cap M|=2$, where $M$
denotes the edge set of $P_2\times T$ between the left half and
the right half of $P_2\times T$. Suppose $|E(C)\cap M|=k>2$. It is
obvious that $k$ even. We suppose that $k=2t$ ($t>1$). Then $C$
has the following form:\\
$P(i_0\rightarrow i_1)\cup v_{i_1}'v_{i_1}''\cup P(i_1\rightarrow
i_2)\cup v_{i_2}''v_{i_2}'\cup P(i_2\rightarrow i_3)\cup
v_{i_3}'v_{i_3}''\cup \ldots \cup P(i_{2t-1}\rightarrow
i_{2t})\cup v_{2t}''v_{2t}'\cup P(i_{2t}\rightarrow i_{2t+1})$,\\
where $P(i_j\rightarrow i_{j+1})$ are paths from vertex $v_{i_j}'$
to vertex $v_{i_{j+1}}'$ in $T_1$ when $j$ ($0\leq j\leq 2t$) are
even, and paths from vertex $v_{i_j}''$ to vertex $v_{i_{j+1}}''$
in $T_2$ when $j$ ($0\leq j\leq 2t-1$) are odd, and
$v_{i_{2t+1}}'=v_{i_0}'$. Hence there exists a cycle of the form
$v_{i_0}'\rightarrow v_{i_1}'\rightarrow v_{i_2}'\rightarrow
\ldots \rightarrow v_{i_{2t-1}}'\rightarrow v_{v_{2t}}'\rightarrow
v_{i_0}'$ in $T_1$. This is a contradiction. Hence the claim
holds.
\par The lemma is immediate from the claim.\\
{\bf Corollary 5} Let $T$ be a tree and $T^e$ be an arbitrary
orientation of $T$. Then the orientation $(P_2\times T)^e$ defined
as above is a Pfaffian orientation of $P_2\times T$ under which
every cycle of even length of $P_2\times T$ is oddly oriented in
$(P_2\times T)^e$.\\
{\bf Proof}\ By Lemma 3, $(P_2\times T)^e$ is a Pfaffian
orientation of $P_2\times T$. Hence, by Lemma 2, every nice cycle
in $(P_2\times T)^e$ is oddly oriented. Then Corollary 5 is
immediate from Lemma 4.
\begin{figure}[htbp]
  \centering
  \scalebox{0.75}{\includegraphics{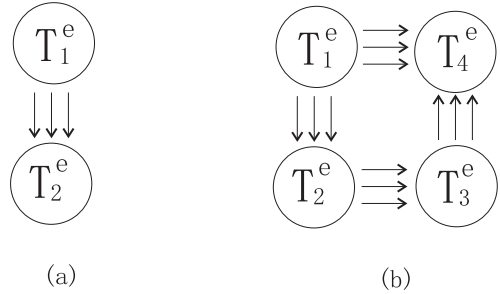}}
  \caption{\ (a)\ An orientation $(P_2\times T)^e$ (=$G^e$).
  \ (b)\ The corresponding orientation $(C_4\times T)^e$. }
\end{figure}
\par Let $T$ be a tree and $T^e$ an arbitrary orientation of $T$. Then, by Corollary 5,
the orientation $(P_2\times T)^e$ defined as above is a Pfaffian
orientation of $P_2\times T$ under which every cycle of even
length of $P_2\times T$ is oddly oriented in $(P_2\times T)^e$.
Let $G=P_2\times T$ and $G^e=(P_2\times T)^e$ and let ${\bf
(P_2\times P_2\times T)^e}=(P_2\times G)^e$ be the orientation
defined as above. Figure 2 illustrates
this procedure, where both of $T_1^e$ and $T_3^e$ are $T^e$,
and both of $T_2^e$ and $T_4^e$ are the converse of $T^e$. Note that $P_2\times P_2=C_4$. Hence
we have $(C_4\times T)^e=(P_2\times P_2\times T)^e$.
From Lemma 3 and Corollary 5, the following theorem is immediate.\\
{\bf Theorem 6} Let $T$ be a tree and let $T^e$ be an arbitrary orientation of $T$.
Then the orientation $(C_4\times T)^e$ of $C_4\times T$ defined as above
is a Pfaffian orientation.\\
{\bf Lemma 7} Let $T$ be a tree, and $T^e$ an arbitrary
orientation. Then $\theta$ is an eigenvalue of $A(T)$ with
multiplicity $m_{\theta}$ if and only if $i\theta$ is an
eigenvalue of $A(T^e)$ with multiplicity $m_{\theta}$, where
$A(T)$ and $A(T^e)$ are the adjacency matrix of $T$
and the skew adjacency matrix of $T^e$, respectively, and $i^2=-1$.\\
{\bf Proof} Let $\phi (T,x)=\det (xI-A(T))$. Since $T$ is a bipartite graph, we may assume
that
\begin{equation}
\phi (T,x)=x^n-a_1x^{n-2}+a_2x^{n-4}+\ldots
+(-1)^ia_ix^{n-2i}+\ldots +(-1)^ra_rx^{n-2r},
\end{equation}
where $n$ and $r$ are the number of vertices of $T$ and the
maximum number of edges in a matching of $T$ (see Biggs [1]). Note
that $(-1)^ia_i$ equals the sum of all principal minors of $A(T)$
of order $2i$. Hence $(-1)^ia_i$ equals the sum of $\det A(H)$
over all induced subgraphs $H$ of $T$ with $2i$ vertices, where
$A(H)$ is the adjacency matrix of subgraph $H$. Note that every
induced subgraph $H$ of $T$ with $2i$ vertices is either a subtree
of $T$ or some subtrees of $T$. Hence $\det A(H)$ equals $(-1)^i$
if $H$ has a perfect matching and 0 otherwise. Thus we have proved
the following claim.\\
{\bf Claim 1} Every $a_i$ equals the number of the induced
subgraphs of $T$ with $2i$ vertices that have a perfect matching.
\par Note that the coefficient of $x^{n-k}$ in $\det (xI-A(T^e))$ is
equal to the sum of $(-1)^k\det A(H^e)$ over all induced
subdigraphs $H^e$ of $T^e$ with $k$ vertices, where $A(H^e)$ is
the skew adjacency matrix of subdigraph $H^e$. It is obvious that
if $k$ is odd then the coefficient of $x^{n-k}$ in $\det
(xI-A(T^e))$ equals 0. Suppose $k$ is even. Let $H$ be the
underlying graph of $H^e$, which is either a subtree of $T$ or
some subtrees of $T$. It is clear that $H^e$ is a Pfaffian
orientation of $H$. Hence $\det A(H^e)$ equals 1 if $H$ has a
perfect matching and 0 otherwise. This implies that the
coefficient of $x^{n-2i}$ in $\det (xI-A(T^e))$ equals the number
of induced subgraphs of $T$ with $2i$ vertices which have a
perfect matching.
Hence we have proved the following claim.\\
{\bf Claim 2} For the orientation $T^e$ of $T$, we have
\begin{equation}
\det (xI-A(T^e))=x^n+a_1x^{n-2}+a_2x^{n-4}+\ldots
+a_ix^{n-2i}+\ldots +a_rx^{n-2r}.
\end{equation}
\par Note that the spectrum of
a bipartite graph is symmetric with respect to $0$ ( see Coulson
and Rushbrooke [3] or Biggs [1]). Hence,  by Claims 1 and 2,
the lemma follows.\\
{\bf Theorem 8} Let $T$ be a tree with $n$ vertices. Then
$$\mbox{Pm}(C_4\times T)=\prod_{j=1}^{n} (2+\theta_j^2),$$
where the eigenvalues of $T$ are $\theta_1, \theta_2, \ldots, \theta_n$.\\
{\bf Proof} Suppose that $(C_4\times T)^e$ is the Pfaffian orientation of
$C_4\times T$ defined as that in Theorem 6. Let $A(T^e)$ be the skew adjacency matrix of $T^e$.
By a suitable labelling of vertices of $(C_4\times T)^e$,
the skew adjacency matrix of $(C_4\times T)^e$
has the following form:
$$A((C_4\times T)^e)=\left [\begin{array}{cccc}
A(T^e) & I & I & 0\\
-I & -A(T^e) & 0 & I\\
-I & 0 & -A(T^e) & -I\\
0 & -I & I & A(T^e)
\end{array}
\right ],
$$
where $I$ is the identity matrix. Hence, by Lemma 1, we have
$$
[\mbox{Pm}(C_4\times T)]^2=\det A((P_2\times P_2\times T)^e)
$$
$$=\det \left [\begin{array}{cccc}
A(T^e) & I & I & 0\\
-I & -A(T^e) & 0 & I\\
-I & 0 & -A(T^e) & -I\\
0 & -I & I & A(T^e)
\end{array}
\right ]
$$
$$
=\det \left \{-\left [\begin{array}{cc}
A(T^e) & I\\
-I & -A(T^e)
\end{array}
\right ]^2+
\left [\begin{array}{cc}
 I & 0\\
0 & I
\end{array}
\right]\right \}
$$
$$
=\det \left [\begin{array}{cc}
2I-(A(T^e))^2 & 0\\
0 & 2I-(A(T^e))^2
\end{array}
\right ].
$$
Hence we have proved
$$\mbox{Pm}(C_4\times T)=|\det (2I-(A(T^e))^2)|.$$
Hence, by Lemma 7, we have
 $$\mbox{Pm}(C_4\times T)=|\det (2I-(A(T^e))^2)|=\prod_{j=1}^n (2+\theta_j^2),$$
where the eigenvalues of $T$ are $\theta_1, \theta_2, \ldots, \theta_n$. The theorem is thus proved.\\
{\bf Remark 9} Note that if $T$ is a path with $n$ vertices, then the set of eigenvalues of $T$ is
$\{2\cos \frac {k\pi}{n+1}| 1\leq k\leq n\} $. Hence, by Theorem 8, the number of perfect matchings of
$C_4\times P_n$ (the linear $2\times 2\times n$ cubic lattice)
equals $\prod\limits_{k=1}^{n} \left [2+4\cos^2 \frac {k\pi}{n+1}\right ]$.
This formula was previously obtained by H. Narumi and H. Hosoya in [14]. \\
{\bf Corollary 10} Suppose $T$ is a tree with $n$ vertices. Then
$\mbox{Pm}(C_4\times T)$ is always a square or double a square.
Moreover, if $T$ is a tree with a perfect matching,
then $\mbox{Pm}(C_4\times T)$ is always a square.\\
{\bf Proof} Suppose that $\phi (T,x)$ is the characteristic polynomial of $T$. Since $T$ is a bipartite graph,
the zeroes of $\phi (T,x)$ are symmetric with
respect to zero
(a result obtained by Coulson and
Rushbrooke [3], see also Biggs [1]). Without loss of generality, we may suppose that
\begin{equation}
\phi (T,x)=x^n-a_1x^{n-2}+a_2x^{n-4}+\ldots
+(-1)^ja_jx^{n-2j}+\ldots +(-1)^ra_rx^{n-2r},
\end{equation}
where $r$ is the number of edges in a maximum matching of $T$. Let
$s=n-2r$. Thus, we have
\begin{equation}
\phi (T,x)=x^s\prod_{j=1}^{r} (x-\theta_j) (x+\theta_j),
\end{equation}
where $\pm \theta_j$ for $1\leq j\leq r$ are all of non$-$zero
eigenvalues of $T$. Hence, by Theorem 8, we have
\begin{equation}
\mbox{Pm}(C_4\times T)=2^s\prod_{j=1}^{r} (2+\theta_j^2)^2.
\end{equation}
Note that $\phi (T, i\sqrt {2})=(i\sqrt {2})^s\prod\limits_{j=1}^{r} (i\sqrt {2}-\theta_j)
(i\sqrt {2}+\theta_j)=(-1)^r(i\sqrt {2})^s\prod\limits_{j=1}^{r} (2+\theta_j^2)$,
where $i^2=-1$.
Hence we have
\begin{equation}
\phi^2(T, i\sqrt 2)=(-1)^s2^s\prod_{j=1}^{r}(2+\theta_j^2)^2.
\end{equation}
By equations (5) and (6), we have
\begin{equation}
\mbox{Pm}(C_4\times T)=2^s\prod_{j=1}^{r}
(2+\theta_j^2)^2=(-1)^s\phi^2 (T, i\sqrt 2).
\end{equation}
Note that $a_j$, for $1\leq j\leq r$, is a non-negative integer, then
$(-1)^s\phi^2(T, i\sqrt 2)
$ equals
$$
(-1)^s\{(i\sqrt 2)^s[(i\sqrt 2)^{2r}-a_1(i\sqrt
2)^{2r-2}+a_2(i\sqrt 2)^{2r-4}+\ldots + (-1)^ja_j(i\sqrt
2)^{2r-2j}+\ldots +(-1)^ra_r]\}^2,
$$
which is a square or double a square. This implies that
$\mbox{Pm}(C_4\times T)$ is a square or double a square. Hence the
first assertion in Corollary 10 holds. If $T$ is a tree with a
perfect matching, then $s=0$ and hence $(-1)^s\phi^2(T, i\sqrt 2)
$ is a square. Thus the second assertion in Corollary 10 holds. The corollary is thus proved.\\
\section* { \bf 3. Enumeration of perfect matchings of $P_3\times T$ and $P_4\times T$ }
\par Suppose that $T$ is a tree and $P_m$ is a path with $m$ vertices.
Let $T^e$ be any orientation of $T$ and let $T_*^e$ be the converse of $T^e$
which is the digraph obtained from $T^e$ by reversing the orientation of each arc.
We define an orientation of $P_m\times T $ (denoted $(P_m\times T)^e$) as follows.
\par Let $V(T)=\{v_1, v_2, \ldots, v_n\}$ be the vertex$-$set of $T$. Take $m$ copies of $T$, denoted by $T_1, T_2,
\ldots, T_m$, where $V(T_i)=\{v_1^{(i)}, v_2^{(i)}, \ldots,
v_n^{(i)}\}$ is the vertex-set of $T_i$, $i=1, 2, \ldots, m$.
Clearly, the mapping $\phi_i$ (from $T$ to $T_i$): $v_j\longmapsto
v_j^{(i)} ( 1\leq j\leq n )$ is an isomorphism between $T$ and
$T_i$. If we add the set of edges $\{v_j^{(i)}v_j^{(i+1)}| 1\leq
j\leq n\}$ between every pair of trees $T_i$ and $T_{i+1}$ for
$1\leq i\leq m-1$, then the resulting graph is $P_m\times T$. We
define the orientation of $T_i$ in $P_m\times T$ to be $T^e$ if
$i$ is odd, denoted by $T_i^e$, and the converse $T_*^e$
otherwise, denoted also by $T_i^e$, and the direction of edges of
the form $v_j^{(i)}v_j^{(i+1)} ( 1\leq j\leq n, 1\leq i\leq m-1 )$
in $P_m\times T$ is from $v_j^{(i)}$ to $v_j^{(i+1)}$. Hence we
obtain an orientation of $P_m\times T$, denoted by $(P_m\times
T)^e$ (see Figure 3).
\begin{figure}[htbp]
  \centering
  \scalebox{0.75}{\includegraphics{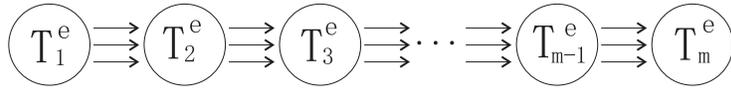}}
  \caption{An orientation $(P_m\times T)^e$ of $P_m\times T$.}
\end{figure}
\par For the sake of convenience, we need introduce some notations.
Let $T$ be a tree with $n$ vertices. For the graph $P_{2k}\times
T$, let $M=M_1\cup M_2\cup \ldots \cup M_k$, where $M_i=
\{v_j^{(2i-1)}v_j^{(2i)}| 1\leq j\leq n \}$ for $1\leq i\leq k$.
Clearly $M$ is a perfect matching of $P_{2k}\times T$. Suppose $T$
is a tree with $n$ vertices containing a perfect matching. For the
graph $P_{2k+1}\times T$, let $M^*=M_1\cup M_2\cup \ldots \cup
M_k\cup M'$, where $M_i= \{v_j^{(2i-1)}v_j^{(2i)}| 1\leq j\leq n
\}$ for $1\leq i\leq k$ and $M'$ is the unique
perfect matching of $T_{2k+1}$. Then $M^*$ is a perfect matching of $P_{2k+1}\times T$. \\
{\bf Lemma 11} Let $T$ be a tree. Then $(P_4\times T)^e$ defined as above is
a Pfaffian orientations of $P_4\times T$. \\
{\bf Proof} Let $M, M_1$ and $M_2$ be defined as above and let $C$ be an $M-$alternating cycle in $P_4\times T$.
By Lemma 2, we only need to prove that $C$ is oddly oriented in $(P_4\times T)^e$. Noting the definitions of
$C_4\times T$ and $P_4\times T$, every nice cycle in $P_4\times T$ is also a nice cycle in $C_4\times T$. Hence
$C$ is a nice cycle in $C_4\times T$.
By the definitions of $(C_4\times T)^e$ and $(P_4\times T)^e$, $(P_4\times T)^e$ is a subdigraph of $(C_4\times T)^e$.
Since $(C_4\times T)^e$ is a Pfaffian orientation of $C_4\times T$, every nice cycle in $C_4\times T$ is oddly oriented
in $(C_4\times T)^e$. Thus $C$ is oddly oriented in $(P_4\times T)^e$. The lemma thus follows.\\
{\bf Lemma 12} Let $T$ be a tree with a perfect matching. Then $(P_3\times T)^e$ defined as above is
a Pfaffian orientations of $P_3\times T$. \\
{\bf Proof} In Lemma 11 we proved that $(P_4\times T)^e$ is a Pfaffian orientation of $P_4\times T$. Note that $T$
contains a perfect matching. Hence every nice cycle in $P_3\times T$ is also a nice cycle in $P_4\times T$. By using
the same method as in Lemma 11, we may prove that $(P_3\times T)^e$ is a Pfaffian orientation of $P_3\times T$.
The lemma is thus proved.\\
{\bf Theorem 13} Suppose $T$ is a tree with $n$ vertices. Then
$$\mbox{Pm}(P_4\times T)=\prod (1+3\alpha^2+\alpha^4),$$ where the
product ranges over all non-negative
eigenvalues $\alpha$ of $T$.\\
{\bf Proof} By Lemma 11, $(P_4\times T)^e$ defined as above is a Pfaffian orientation of $P_4\times T$. Hence,
by Lemma 1, we have
$$[\mbox{Pm}(P_4\times T)]^2=\det [A((P_4\times T)^e)],$$
where $A((P_4\times T)^e)$ is the skew adjacency matrix of $(P_4\times T)^e$.
By a suitable labelling of vertices of $(P_4\times T)^e$,
the skew adjacency matrix of $(P_4\times T)^e$
has the following form:
$$A((P_4\times T)^e)=\left [\begin{array}{cccc}
A & I & 0 & 0  \\
-I & -A & I &  0   \\
0 & -I & A & I     \\
0 & 0 & -I & -A   \\
\end{array}
\right ],
$$
where $A$ denotes the skew adjacency matrix $A(T^e)$ of $T^e$.
\par Now multiplying the first column, then the third and fourth row, then the fourth column
of the partitioned matrix $A((P_4\times T)^e)$ by $-1$, we do not
change the absolute value of the determinant and we obtain matrix
$Q$, where
$$Q=\left [\begin{array}{cccc}
-A & I & 0 & 0  \\
I & -A & I & 0 \\
0 & I & -A & I   \\
0 & 0 & I & -A   \\
\end{array}
\right ].
$$
Denote by $B$ the adjacency matrix of the path with four vertices,
that is,
$$B=\left [\begin{array}{cccc}
0 & 1 & 0 & 0 \\
1 & 0 & 1 & 0   \\
0 & 1 & 0 & 1 \\
0 & 0 & 1 & 0 \\
\end{array}
\right ].
$$
Then we may write
$$Q=-I_4\otimes A+B\otimes I_n,$$
where $\otimes$ denotes the Kronecker product of matrices.
\par Note that, if $A$ has the eigenvalues
$\lambda_1, \lambda_2, \ldots, \lambda_n$ and $B$ has the
eigenvalues $\mu_1, \mu_2, \mu_3$ and $\mu_4$, then the
eigenvalues of $-I_4\otimes A+B\otimes I_n$ are as follows:
\begin{center}
$\mu_i-\lambda_j$, \ \ \ \mbox{where} \ \ \ $1\leq i\leq 4, \ \ \
1\leq j\leq n$.
\end{center}
\par Suppose that $T$ has the eigenvalues $\alpha_1, \alpha_2, \ldots, \alpha_n$.
By Lemma 7, $A$ has the eigenvalues $i\alpha_j$ ($1\leq j\leq n$
), where $i^2=-1$. Note that the eigenvalues of $B$ are as
follows:
$$\pm \sqrt {\frac {3+\sqrt {5}}{2}},\ \ \pm \sqrt {\frac {3-\sqrt {5}}{2}}.$$
Thus the eigenvalues of $Q$ are as follows:
\begin{center}
$\pm \sqrt {\frac {3+\sqrt {5}}{2}}-i\alpha_s, \pm \sqrt {\frac
{3-\sqrt {5}}{2}}-i\alpha_s, (s=1, 2, \ldots, n).$
\end{center}
Hence the determinant of matrix $Q$ is the product of these numbers. Since we are interested in the absolute
value of this determinant, we may replace these $4n$ factors by their absolute values, and so the
absolute value of the determinant of matrix $A((P_4\times T)^e)$ is
$$\left |\prod_{s=1}^{n} \left (\sqrt {\frac {3+\sqrt {5}}{2}}-i\alpha_s\right )
\left (-\sqrt {\frac {3+\sqrt {5}}{2}}-i\alpha_s\right) \left
(\sqrt {\frac {3-\sqrt {5}}{2}}-i\alpha_s\right )\left (-\sqrt
{\frac {3-\sqrt {5}}{2}}-i\alpha_s\right) \right |
$$
$$
=\prod_{s=1}^{n} (1+3\alpha_s^2+\alpha_s^4).
$$
Hence
$$\mbox{Pm}(P_4\times T)=
\prod_{s=1}^{n} (1+3\alpha_s^2+\alpha_s^4)^{\frac {1}{2}}.
$$
Note that the spectrum of a tree is symmetric with respect to zero
( see Coulson and Rushbrooke [3] or Biggs [1]). Hence we have
$$\mbox{Pm}(P_4\times T)=
\prod (1+3\alpha^2+\alpha^4),
$$
where the product ranges over all non-negative eigenvalues of $T$.
The theorem is thus proved.
\par Similarly, by using Lemma 12, we may prove the following theorem.\\
{\bf Theorem 14} Suppose $T$ is a tree with a perfect matching. Then
$$\mbox{Pm}(P_3\times T)=\prod (2+\alpha^2),$$
where the product ranges over all positive eigenvalues $\alpha$ of $T$.\\
{\bf Corollary 15} Suppose $T$ is a tree with a perfect matching.
Then $[\mbox{Pm}(P_3\times T)]^2=\mbox{Pm}(C_4\times T)$.
\par Corollary 15 is immediate from Theorems 8 and 14.
\par Although a tree $T$ with even number of vertices has no perfect matching, $P_3\times T$ may
contain perfect matchings. See for example the tree $T$ in Figure
4, which has no perfect matching but $P_3\times T$ contains a
perfect matching (the set of the bold edges).
\begin{figure}[htbp]
  \centering
  \scalebox{0.75}{\includegraphics{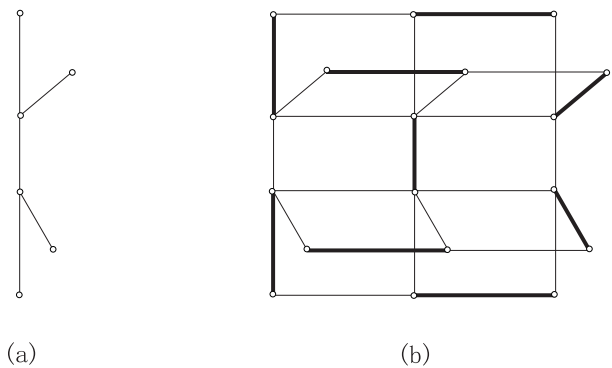}}
  \caption{\ (a)\ A tree $T$ having no perfect matching.
  \ (b)\ $P_3\times T$.}
\end{figure}
Hence we pose naturally the following problems.\\
{\bf Problem 1} Suppose that $T$ is a tree with even number of vertices containing no perfect matching. Enumerate
perfect matchings of $P_3\times T$.\\
{\bf Problem 2} Suppose that $T$ is a tree and $m>4$. Enumerate perfect matchings of $P_m\times T$.\\
{\bf Remark 17} If the tree in Problem 2 is a path $P_n$, then the number of perfect matchings of $P_m\times P_n$ equals
$$2^{\frac {mn}{2}}\prod_{k=1}^{m}\prod_{l=1}^{n} \left ( \cos ^2\left (\frac {\pi k}{m+1}\right )+
\cos ^2\left (\frac {\pi l}{n+1}\right )\right )^{\frac {1}{4}},$$
which was obtained by a physicist, Kasteleyn (see [8,9,11]). It is well known as the dimer problem,
which has applications in statistical mechanics.\\
\vskip0.5cm
\noindent
{\bf Acknowledgements}
\par We wish to thank Professor Richard Kenyon for some useful discussions.\\
\vskip0.5cm
\newcounter{cankao}
\begin{list}
{[\arabic{cankao}]}{\usecounter{cankao}\itemsep=0cm}
\centerline{\bf References}
\vspace*{0.5cm}
\small
\item N. Biggs, Algebraic Graph Theory, Cambridge, Cambridge University Press, 1993.
\item M. Ciucu, Enumeration of perfect matchings in graphs with reflective symmetry,
J. Combin. Theory Ser. A 77(1997), 67$-$97.
\item C. A. Coulson and G. S. Rushbrooke, Note on the method of molecular orbits, Proc. Camb. Philos. Soc. 36(1940),
193$-$200.
\item S. J. Cyvin and I. Gutman, Kekul\' e structures in Benzennoid Hydrocarbons, Springer
Berlin, 1988.
\item I. Fischer and C. H. C. Little, Even Circuits of Prescribed Clockwise Parity, Preprint.
\item G. G. Hall, A Graphic Model of a Class of Molecules, Int. J. Math. Edu. Sci. Technol., 4(1973), 233$-$240.
\item W. Jockusch, Perfect matchings and perfect squares,
J. Combin. Theory Ser. A 67(1994), 100$-$115.
\item P. W. Kasteleyn, Dimer statistics and phase transition, J. Math. Phys. 4(1963), 287$-$293.
\item P. W. Kasteleyn, Graph Theory and Crystal Physics. In F.Harary, editor, Graph
Theory and Theoretical Physics. Academic Press, 1967, 43$-$110.
\item C. H. C. Little, A characterization of convertible (0,1)$-$matrices,
J. Combinatorial Theory 18(1975), 187$-$208.
\item L. Lov\' asz and M. Plummer, Matching Theory, Ann. of Discrete Math. 29, North$-$Holland,
New York, 1986.
\item W. McCuaig, P\' olya's permanent problem, Preprint.
\item W. McCuaig, N. Robertson, P. D. Seymour, and R. Thomas, Permanents, Pfaffian orientations,
and even directed circuits (Extended abstract), Proc. 1997 Symposium on the
Theory of Computing (STOC).
\item H. Narumi and H. Hosoya, Proof of the generalized expressions of the number of perfect matchings of polycube
graphs, J. Math. Chem. 3(1989), 383$-$391.
\item L. Pauling, The Nature of Chemical Bond, Cornell. Univ. Press, Ithaca, New York, 1939.
\item J. Propp, Enumeration of Matchings: Problems and Progress, In: New Perspectives in Geometric
Combinatorics (eds. L. Billera, A. Bj\"orner, C. Greene, R.
Simeon, and R. P. Stanley), Cambridge University Press, Cambridge,
(1999), 255$-$291.
\item N. Robertson, P. D. Seymour, and R. Thomas, Permanents, Pfaffian orientations, and even
directed circuits, Annals of Math., 150(1999), 929$-$975.
\item R. Swinborne$-$Sheldrake, W. C. Herndon and I. Gutman, Kekul\' e structures and resonance energies of
benzennoid hydrocarbons, Tetrahedron Letters, (1975), 755$-$758.
\item W. Yan and F. Zhang, Enumeration of perfect matchings of graphs with reflective symmetry by Pfaffians,
Adv. Appl. Math., 32(2004), 655$-$668.
\item W. Yan and F. Zhang, On the Number of Kekul\' e Structures of a Type of Oblate Rectangles,
MATCH$-$ Commun. Math. Comput. Chem., 47(2003), 141$-$149.
\item F. Zhang and W. Yan, Enumeration of perfect matchings in a type of graphs with reflective symmetry,
MATCH$-$ Commun. Math. Comput. Chem., 48(2003), 117$-$124.
\item F. Zhang and H. Zhang, A new enumeration method for Kekul\'e structures of hexagonal systems with
forcing edges, J.Mol.Struct. (THEOCHEM) 331(1995), 255$-$260.
\item F. Zhang and H. Zhang, A note on the number of perfect matchings of bipartite graphs, Discrete Appl. Math.
73 (1997), 275$-$282.
\end{list}
\end{document}